\documentclass[11pt]{article}
\usepackage{amsmath,amsthm,amssymb}
\usepackage{amsfonts}
\usepackage{graphics}
\usepackage{color}
\newtheorem{intro}{Theorem}[section]

\newtheorem{theorem}{Theorem}[subsection]
\newtheorem{lemma}[theorem]{Lemma}

\newtheorem{proposition}[theorem]{Proposition}
\theoremstyle{definition}

\newtheorem{remark}[theorem]{Remark}

\newtheorem{example}[theorem]{Example}

\def\pri{{\mathcal{P}}{\mathbb{Z}}}
\def\slnr{SL(n,{\mathbb{R}})}   
\def\slsr{SL(s,{\mathbb{R}})}   
\def\slnz{SL(n,{\mathbb{Z}})}   
\def\slsz{SL(s,{\mathbb{Z}})}   
\def\ma{\mathrm{max}}
\def\calp{\mathcal{P}}   
\def\calr{\mathcal{R}}   
\def\tt{\mathcal{T}}   
\def\tts{\tt_s}   
\def\pp{\mathcal{P}}   
\def\pps{\pp_s}   
\def\mm{\mathcal{M}}   
\def\hh{\mathcal{H}}   
\def\cl{\mathcal{L}}   
\def\tl{\widetilde{\mathcal{L}}}   
\def\tm{\widetilde{\mathcal{M}}}   
\def\bu{\mathbf{u}}
\def\1n{\frac{1}{n}} 

\def\bf1{\bar{f}_1}
\def\bfs{\bar{f}_{s-1}}

\newcommand {\iv}{^{-1}}
\newcommand {\dist}{\mathrm{dist}}
\newcommand {\pr}{\mathrm{proj}}
\newcommand {\N}{\mathbb{N}} 
\newcommand {\Z}{\mathbb{Z}}            
\newcommand {\R}{\mathbb{R}} 
\newcommand {\Q}{\mathbb{Q}} 
\newcommand {\C}{\mathbb{C}} 
\newcommand {\K}{\mathbb{K}} 
\newcommand {\g}{\frak g} 

\newcommand {\is}{i\in \{ 1,2,\dots s-1 \}}
\newcommand {\me}{\medskip}

\newcommand {\bi}{\bigskip}

\newcommand {\Notat}{\noindent {\it{Notation}}:\thickspace } 
\begin{document}
\makeatletter
\title{Transference principles and locally symmetric spaces}
\author{
Cornelia DRU\c{T}U\\ \\
Mathematical Institute,\\
24-29 St Giles,\\
Oxford, OX1 3LB\\
United Kingdom\\
drutu@maths.ox.ac.uk}
\date{ }
\maketitle

\begin{abstract}
\noindent We explain how the Transference Principles from Diophantine approximation can be interpreted in terms of geometry of the
locally symmetric spaces ${\mathcal{T}}_n=SO(n) \backslash SL(n,{\mathbb{R}}) /SL(n,{\mathbb{Z}})$ with  $n\geq 2$, and how, via this dictionary, they become transparent geometric remarks and can be easily proved. Indeed, a finite family of linear forms is naturally identified to a locally geodesic ray in a space ${\mathcal{T}}_n$ and the way this family is approximated is reflected by the heights at which the ray rises in the cuspidal end. The only difference between the two types of approximation appearing in a Transference Theorem  is that the height is measured with
respect to different rays in $\overline{W}_0$, a Weyl chamber in ${\mathcal{T}}_n$. Thus the Transference Theorem is equivalent to a relation between the Busemann functions of two rays in $\overline{W}_0$. This relation is easy to establish on $\overline{W}_0$, because restricted
to it the two Busemann functions become two linear forms. Since ${\mathcal{T}}_n$ is at finite Hausdorff
distance from $\overline{W}_0$, the same relation is satisfied up to a bounded perturbation on the whole of
${\mathcal{T}}_n$.
\end{abstract}

\tableofcontents

\section{Introduction}

In this paper we explain how the Transference principles appearing in Diophantine approximation of
systems of linear forms have an easy interpretation in terms of geometry of the locally symmetric
spaces $SO(n)\backslash \slnr / \slnz$.

Consider a family of $\ell$ linear real forms in $m$ variables,  and the transposed family of $m$
linear forms in $\ell$ variables:
$$
L_i(x_1,...,x_m) =\sum_{j=1}^m a_{ij} x_j\, ,\; \;  M_j(y_1,..., y_\ell) = \sum_{i=1}^\ell a_{ij} y_i
\, .
$$

We denote by $L$ the $\ell\times m$ matrix $L = (a_{ij})_{1\le i\le\ell,_, 1\le j\le m}$ and by $M$ its transpose.

We also denote throughout the paper by $\|\cdot \|_e$ the Euclidean norm and by $\|\cdot \|_\ma$ the
max-norm in $\R^n$, that is the norm defined by:
$$\|x \|_\ma = \max \{|x_1|, |x_2|, \ldots , |x_n| \}\, .$$

We denote by $\pri^n $ {\it{the set of primitive integer vectors in }}$\R^n$,
$$\{ (p_1,p_2,\dots ,p_n)\in
\Z^n \setminus \{ (0,\dots 0) \} \; ;\;  {\rm{gcd}}(p_1,p_2,\dots ,p_n)=1 \}\, .$$



Dirichlet's theorem implies the existence of infinitely
many integral solutions 
$$(\bar{p}, \bar{q})= (p_1,..., p_\ell, q_1,..., q_m )\, $$ for the following equation:
\begin{equation}\label{dir}
\| L(\bar{q}) - \bar{p} \|_\ma \leq  \| \bar{q} \|_\ma^{-m/\ell}\, .
\end{equation}

The family of linear forms $(L_i)_{1\leq i \leq \ell }$ is said to be \emph{very well approximable} if for some
$\alpha >0$ and infinitely many integral vectors $(\bar{p}, \bar{q})$ the following holds:
\begin{equation}\label{vwa}
\| L(\bar{q}) - \bar{p} \|_\ma \leq  \| \bar{q} \|_\ma^{-(m+\alpha)/\ell }\, .
\end{equation}

The definition of very well approximable forms can be slightly generalized using approximating
functions. Throughout, an \emph{approximating function} is a decreasing function $\phi :\R_+ \to \R_+$
with $\lim_{x\to \infty} \phi (x) =0$. In (\ref{vwa}) one can replace the second term by $\phi (\|
\bar{q} \|_\ma )$, where $\phi$ is an approximating function such that $\lim_{x\to \infty }
x^{\frac{\ell }{m }}\phi (x) =0$.

The transference principles state that if $(L_i)_{1\leq i \leq \ell }$ is $\phi$-very well approximable
then the transposed family $(M_j)_{1\leq j \leq m }$ is $\psi$-very well approximable, and give
estimates of $\psi$ in terms of $\phi$. More precisely, the following theorem is a version in terms
 of approximating functions of \cite[Theorem II, Chapter V,
Section 2]{Cassels:diophantine}:

\begin{intro}\label{functgenkhintch}
\begin{itemize}
    \item[(I)] Assume that the following system of inequalities has infinitely many integral solutions $(\bar{p},
\bar{q})= (p_1,..., p_\ell, q_1,..., q_m )$:
\begin{equation}\label{l}
|L_i(\bar{q}) - p_i |\leq \phi \left( \| \bar{q} \|_\ma \right) \, .
\end{equation}

Then there exist infinitely many integral solutions $(\bar{a}, \bar{b})= (a_1,..., a_\ell, b_1,..., b_m
)$ for the system of inequalities:
\begin{equation}\label{m}
|M_j(\bar{a}) - b_i |\leq \psi \left( \| \bar{a} \|_\ma \right)\, ,
\end{equation}
where $\psi = F \circ G\iv$ with $F(x)=(s-1)x^{\frac{1-\ell}{s-1}}\phi(x)^{\frac{\ell}{s-1}}$ and
$G(x)=(s-1)x^{\frac{m}{s-1}}\phi(x)^{\frac{1-m}{s-1}}$, and $s=\ell +m$.

\me

    \item[(II)] In particular, if $\phi (x) = x^{-\frac{m+\alpha }{\ell}}$ then $\psi (x) =  x^{-\frac{\ell+\beta
}{m}}$ with $\beta = \frac{\ell \alpha }{m(m+\ell -1)+ (m-1)\alpha }\, .$
\end{itemize}
\end{intro}

Note that if $\phi$ is an approximating function then the function $F:\R_+ \to \R_+$ is decreasing and
$\lim_{x\to \infty } F(x)=0$ while the function $G:\R_+ \to \R_+$ is increasing and $\lim_{x\to \infty
} G(x)=+\infty$. The two properties imply that $\psi$ is an approximating function.

By applying Theorem \ref{functgenkhintch} first for $\ell =1 \, ,\, m=n$ and then for $m=1\, , \, \ell
=n$ the following well-known result is obtained:

\begin{intro}[Khintchine's transference principle]
Let $v=(x_1,...,x_m)$ be a vector in $\R^n$ with all coordinates irrational, and let $\omega (v)$ and
$\omega'(v)$ be the respective least upper bounds of the real numbers $\alpha$ and $\alpha'$ for which
the following inequalities have infinitely many integer solutions $(p,q_1,...,q_n)$ and respectively
$(p_1,...,p_n,q)$:
\begin{equation}\label{khintch}
\begin{array}{cc}
  |q_1x_1+\cdots q_n x_n -p | & \leq  \| (q_1,...,q_n)\|_\ma^{-n-\alpha}\, ,\\
  \max_{1\leq j\leq n} |qx_j-p_j| & \leq q^{-\frac{1+\alpha'}{n}}\, . \\
\end{array}
\end{equation}
Then the following sequence of inequalities holds (with $\omega (v)$ and/or $\omega'(v)$ possibly
taking the value $+\infty $):
\begin{equation}\label{khintchin}
 \frac{\omega (v)}{n^2+(n-1)\omega (v)}\leq \omega'(v) \leq \omega (v)\, .
\end{equation}
\end{intro}

\me

The first to notice a relationship between Diophantine approximation of systems of linear forms and
behavior of geodesic rays in locally symmetric spaces was Dani \cite{Dani:divergent}. He noticed that
systems of linear forms can be identified with unipotent elements in some group $\slnr$, and that the
way in which the systems of linear forms are approximated reflects the way in which a locally geodesic
ray naturally associated to the unipotent travels in the cusp. We shall follow this initial idea to
point out that the Transference theorems become, via this dictionary, very transparent geometric
remarks.

\me

More precisely, let $L\in M_{\ell \times m} (\R )$ and $M\in M_{m\times \ell }(\R )$, and consider the
semisimple group $\slsr$ with $s=\ell +m$, the symmetric space associated to it $\pps = SO(s)
\backslash \slsr$, and the locally symmetric space $\tts = \pps /\slsz$. It is well known that $\tts$
is at finite Hausdorff distance from the isometric image of a Weyl chamber $\overline{W}_0$, which is
an Euclidean sector of dimension $n-1$ and of shape prescribed by the Dynkin diagram of $\slsr$.

Both matrices $L$ in the space $ M_{\ell \times m} (\R )$, and matrices $M$ in $M_{m\times \ell }(\R )$
can be identified to unipotents in $\slsr$ (see (\ref{lmap}) and (\ref{mmap})), so they can also be
naturally identified to locally geodesic rays in $\tts$. Moreover, when $M=L^T$ the unipotent is the
same, and so is the ray. It is shown 
that $L$ (respectively $M$) is very well approximable if and only if that locally geodesic
ray goes infinitely many times in the cusp at larger and larger heights, with lower bounds on the
height given by a function of the time when the height is attained. The function depends on the
approximation function initially considered. See Proposition \ref{lr} for the precise statements.

The only difference between the case of $L$ and the case of $M$ is that the height is measured with
respect to different rays in $\overline{W}_0$. Indeed $\overline{W}_0$ contains a ray $\bar{r}_1$,
which is equally a 1-dimensional face of $\overline{W}_0$, and whose lifts in $\pps$ have as boundary
at infinity all the rational points. (Recall that the boundary at infinity of $\pps$ can be identified
to the spherical building of flags in $\mathbb{P}^{n-1}\R$.) Another 1-dimensional face of
$\overline{W}_0$, the ray $\bar{r}_{s-1}$, has the property that its lifts have as boundaries at
infinity all the rational hyperplanes in $\mathbb{P}^{n-1}\R$.

When considering $L$, the height of the ray in $\tts$ has to be measured with respect to $\bar{r}_1$,
that is, using the Busemann function of  $\bar{r}_1$ (see Section 2.1 for a definition). When studying
$M$, the height must be measured with respect to $\bar{r}_{s-1}$. A picture representing the case $s=3$
can be seen in Figure \ref{rphi}.

It follows that in order to relate an approximating function for $L$ to an approximating function for
$M=L^T$ one needs to relate the Busemann function of $\bar{r}_1$ to the  Busemann function of
$\bar{r}_{s-1}$, on $\tts$. This relation is easy to establish on $\overline{W}_0$, because restricted
there the two Busemann functions become two linear forms (see Example \ref{brn}), and the angle between
their respective vectors of coefficients is completely determined. Since $\tts$ is at finite Hausdorff
distance from $\overline{W}_0$, the same relation up to a bounded perturbation holds on the whole of
$\tts$ (see Proposition \ref{geomtransf} and Figure \ref{quotient}).

This easy to see geometric relation between the two Busemann functions turns out to be the same as the
relation between approximating functions in transference principles. This is illustrated by some
explicit computations in Section \ref{calc}.

\me

The plan of the paper is as follows. In Section \ref{Pr} notations  are introduced, some general
notions and results are recalled, and the objects and formulas from the theory of symmetric spaces are
made explicit in the case of the symmetric space $\pps$. The same is done for the locally symmetric
space $\tts$, moreover in Proposition \ref{geomtransf} an estimate relating the Busemann
functions of $\bar{r}_1$ and $\bar{r}_{s-1}$ is proved.

In Section \ref{sdioph} are described two families of geodesic rays rising in the cusp, it is explained
that their definition is natural, and Proposition \ref{geomtransf} is used to relate them (see Lemma
\ref{relr}).

The relation between sets of very well approximable linear forms and sets of geodesic rays rising in
the cusp is established in Proposition \ref{lr}.

Summing up Proposition \ref{lr} and Lemma \ref{relr} one is able to reprove the transference results.
This is shown in Section \ref{calc} by an explicit computation.

\section{Preliminaries on (locally) symmetric
spaces}\label{Pr}

\subsection{Notations}\label{not}

We denote by diag$(a_1,a_2,\dots ,a_n)$ the diagonal matrix in $\slnr$ having entries $a_1,a_2,\dots
,a_n$ on the diagonal. We denote by $Id_n$ the identity matrix.

We sometimes call a 1-dimensional linear (sub)space of $\R^n$ a \emph{line}.


We denote by $\langle v_1,\dots ,v_k \rangle$ the linear subspace generated by the vectors $v_1,\dots
,v_k$.

\me

Given two functions $f$ and $g$ with values in $\R$, we write $f\ll g$ if $f(x)\leq C\cdot g(x)$, for
every $x$, where $C >0$ is a universal constant. We write $f\asymp g$ if both $f\ll g$ and $f\gg g$
hold.

 We also use the notation $f+O(1)$ to signify a function of the form $f+C$ with $C$ a fixed constant.

\me

If $G$ is a group, we denote by $Z(G)$ its center $\{ z\in G \; ;\; zg=gz\, ,\, \forall g\in G\}$. If
$H$ is a subgroup of $G$ we denote by $C_G (H)$ the center of $H$ in $G$, that is the group $\{ z\in G
\; ;\;  zh=hz\, ,\, \forall h\in H\}$.

If $G$ is a topological group, we denote by $G_e$ its connected component containing the identity.

\me

Let $X$ be a complete Riemannian manifold of non-positive curvature. Two geodesic rays $r$ and $r'$ in
$X$ are called {\it{asymptotic}} (and we write $r\sim r'$) if they are at finite Hausdorff distance one
from the other. The boundary at infinity of $X$ is the quotient $\mathcal{R}/\sim $ of the set
$\mathcal{R}$ of all geodesic rays in $X$ by the equivalence relation $\sim $. It is usually denoted by
$\partial_\infty X$. Given $\xi\in
\partial_\infty X$, we signify that a geodesic ray $r$ is in the equivalence class $\xi$ by the equality $r(\infty )=\xi$.

Let $r$ be a geodesic ray in $X$. {\it The Busemann function associated to $r$} is the function
$$
f_r:X\to \R \, ,\; f_r(x)=\lim_{t\to \infty}[\dist (x,r(t))-t]\; .
$$
The limit exists because the function $t\mapsto \dist (x,r(t))-t$ is non-increasing and bounded.

\begin{lemma}[\cite{BGS}]\label{blip}
For any geodesic ray $r$ in $X$ and any two points $x,y$ in $X$,
$$|f_r(x) - f_r(y)|\leq \dist (x,y)\, .$$
\end{lemma}

The level hypersurfaces $H_a(r)=\lbrace x\in X \; ;\;  f_r(x)= a \rbrace$ are called {\it horospheres},
and the sublevel sets $Hb_a(r)=\lbrace x\in X \; ;\;  f_r(x)\leq a \rbrace$ are called {\it{closed
horoballs}}. For $a=0$ we use the notation $H(r)$ for the horosphere, and $Hb(r)$ for the closed
horoball.

\begin{example}\label{brn}
If $X=\R^n$ an arbitrary geodesic ray is of the form $r(t)= t\cdot v +w $, where $v$ and $w$ are
vectors and $\| v\|_e=1$. An easy computation shows that $f_r(x) = - \langle x|v \rangle + \langle x|w
\rangle$, where $\langle \cdot | \cdot \rangle$ is the standard inner product. In particular
$Hb(r)= \{x\in \R^n \; ;\; \langle x|v \rangle \geq \langle x|w \rangle \}$.
\end{example}

\me

Assume now that the manifold $X$ is also simply connected. The Busemann functions of two asymptotic
rays in $X$ differ by a constant \cite{BridsonHaefliger}. Therefore they are also called {\it{Busemann
functions of basepoint }}$\xi$, where $\xi $ is the equivalence class containing the two rays. The
families of horoballs and horospheres are the same for the two rays. We shall say that they are
horoballs and horospheres {\it of basepoint }$\xi$.


\subsection{The symmetric space $\pps$ of positive definite quadratic forms of determinant one}\label{str}

Throughout the paper we shall identify a quadratic form $Q$ on $\R^s $ with its matrix $M_Q$ in the
canonical basis of $\R^s$. We shall denote by $b_Q$ the bilinear form associated to $Q$.

In what follows we freely use the terminology and the results from the theory of symmetric spaces of
non-compact type without Euclidean factors, and associated semisimple groups. We refer the reader to
\cite{Helgason:diffgeom},  \cite[Chapter 3]{CheegerEbin}, \cite{OnishchikVinberg},
\cite{Mostow:rigidity}, \cite{Raghunathan:discrete} and \cite{Witte:Arithmetic} for details on the
theory.

We study here mainly one such space, that is \emph{the space ${\mathcal{P}}_s$ of positive definite
quadratic forms of determinant one on }$\R^s$. It can be endowed with a metric defined as follows.
Given $Q_1\, ,\, Q_2 \in {\mathcal{P}}_s$, there exists an orthonormal basis with respect to $Q_1$ in
which $Q_2$ becomes diagonal with coefficients $\lambda_1,\dots ,\lambda_s\in \R_+^*$. We define
\begin{equation}\label{dist}
d(Q_1,Q_2)=\sqrt{\sum_{i=1}^s (\ln \lambda_i)^2}\; .
\end{equation}

The connected component of the identity of the group of isometries of ${\mathcal{P}}_s$ can be
identified to the semisimple group $PSL(s, \R )$. This group acts on the right on $\mathcal{P}_s$ by
\begin{displaymath}
\Phi : PSL(s,\R)\times \mathcal{Q}_s\to \mathcal{Q}_s\, ,\, \Phi (B,M)=B^TMB\; .
\end{displaymath} The action can be written in terms of
quadratic forms as $\Phi (B,Q)=Q[B]=Q\circ B$.

Thus, the symmetric space ${\mathcal{P}}_s$ can be identified with
 $SO(s) \backslash SL(s, \R )$ by associating to each
right coset $SO(s)\, Y$ the quadratic form $Q_Y$ whose matrix in the canonical basis is $M_Y=Y^T\cdot
Y$.

The Lie algebra of $SL(s, \R)$ is $\g_s= \{ B\in L(s, \R) \; ;\; \mathrm{trace}\, B =0 \}$.

In the symmetric space $\mathcal{P}_s$ consider as a fixed basepoint the quadratic form $Q_0$ of matrix
$Id_s$. The geodesic symmetry with respect to this point is a global isometry, as $\mathcal{P}_s$ is a
symmetric space. In terms of matrices in the canonical basis of $\R^s$ the symmetry with respect to
$Q_0$ can be written as $M_Q\mapsto M_Q\iv$. It defines on $SL(s, \R)$ the involution $B\mapsto
(B^T)\iv$. The corresponding Cartan involution on the Lie algebra $\g_s$ is $\theta (B)= -B^T$, and the
Cartan decomposition is
 $\g = {\mathfrak k}\oplus {\mathfrak p}$ with ${\mathfrak k}= \{ B\in \g_s \; ;\; B^T=-B
\}$ and ${\mathfrak p}= \{ B\in \g_s \; ;\; B^T=B \}$.

The Killing form on the Lie algebra $\g_s$ is ${\mathfrak Q} (A,B)= \mathrm{trace} (AB)$, hence
$-{\mathfrak Q} (A,\theta (B))= \mathrm{trace} (AB^T)$ defines a positive definite quadratic form on
$\g$ invariant under the adjoint representation $Ad$ restricted to $SO(s)$. The projection $SL(s, \R)
\to {\mathcal{P}}_s = SO(s) \backslash SL(s, \R )$ is a Riemannian submersion.


An element $g_0$ in $\slsr$ is \textit{hyperbolic} if there exists $g\in GL(n,\R)$ such that $g g_0
g\iv$ is diagonalizable with all the eigenvalues real positive.

Consider a field $\K \subset \R$. We say that a Lie group $G$ is \textit{defined over $\K$} if $G$ has
finitely many connected components and if its connected component of the identity coincides with the
connected component of the identity of a real algebraic group defined over $\K$ \cite[Definition
6.2]{Witte:Arithmetic}.

A \textit{torus} is a closed connected Lie subgroup of $\slsr$ which is diagonalizable over $\C$, i.e.
such that there exists $g\in GL(s,\C)$ with the property that $g\, T\, g\iv $ is diagonal. A torus is
called $\K$-\textit{split} if it is defined over $\K$ and diagonalizable over $\K$, that is if there
exists $g\in GL(s,\K)$ with the property that $g\, T\, g\iv $ is diagonal.

\me

\noindent \textit{Conventions}: In this paper by \textit{torus} we mean a non-trivial $\R$-split torus.
By wall/Weyl chamber we mean a \textit{closed} wall/Weyl chamber.
 By its \textit{relative interior} we mean the open wall/Weyl chamber.

\me

 We call {\it{singular torus in}} $\slsr$ a torus $A_0$ which,
 in every maximal torus $A$ containing it,
  can be written as $\bigcap_{\lambda \in \Lambda}\ker \lambda$,
   where $\Lambda $ is a non-empty set of roots on $A$.
    Any such torus is a union of walls of Weyl chambers.

\me

The subgroup of $\slsr$ $$A=\{ \mbox{diag}(e^{t_1},e^{t_2},\dots , e^{t_s})\; ;\;  t_1+t_2+\cdots +
t_s=0 \}$$ is a maximal $\Q $-split torus as well as a maximal $\R $-split torus. A $\Q $-Weyl chamber
(as well as an $\R $-Weyl chamber) is $\triangleleft A= \{ \mbox{diag}(e^{t_1},e^{t_2},\dots ,
e^{t_s})\; ;\; t_1+t_2+\cdots +t_s=0,\; t_1\geq t_2\geq\cdots \geq t_s \}$.

We recall that a \textit{flat} in $\pps$ is a totally geodesically embedded copy of an Euclidean space
in $X$, and that a \textit{maximal flat} is a flat which is maximal with respect to the inclusion.

For instance, the set of positive definite quadratic forms $$F_0=\{ \mbox{diag }(e^{t_1},e^{t_2},\dots
, e^{t_s})\; ;\; t_1+t_2+\cdots + t_s=0 \}$$ is a maximal flat. Note that $F_0$ is nothing else than
the orbit $Q_0[A]$. Finitely many hyperplanes in $F_0$ appear as intersections of it with other maximal
flats through $Q_0$. These hyperplanes split $F_0$ into finitely many Weyl chambers. One of them is the
Weyl chamber $W_0=Q_0[\triangleleft A]$, i.e. the subset of quadratic forms whose matrices moreover
satisfy $t_1\geq t_2\geq\cdots \geq t_s$. The others can be obtained by performing all the possible
permutations in the sequence of inequalities defining $W_0$.

The group $\slsr$ acts transitively on the collection of maximal flats, as well as on the collection of
Weyl chambers in $X$. This is equivalent to saying that it acts transitively by conjugation on the
collection of maximal tori and on the collection of Weyl chambers in $G$.

\me

The dimension 1 walls (singular rays) of $W_0$, parameterized with respect to the arc length, are the
sets of quadratic forms
\begin{equation}\label{rays}
r_i= \{ \mbox{diag }(\underbrace{e^{\lambda_it},\dots ,e^{\lambda_it}}_{s-i\mbox{ times}},
\underbrace{e^{-\mu_i t},\dots  e^{-\mu_i t}}_{i\mbox{ times}})\; ;\;  t\in \R_+ \}\; ,
\end{equation} where $\lambda_i=\sqrt{\frac{i}{s(s-i)}}$ and $\mu_i=\sqrt{\frac{s-i}{si}}\, $,
$i\in \{ 1,2,\dots s-1 \}$.

\me

\subsection{Parabolic and unipotent subgroups of $\pps$}\label{ssg}

There are two ways of defining parabolic subgroups, we recall them both.

If $\triangleleft A_0$ is a wall or a Weyl chamber in the torus $A_0$, the parabolic group
corresponding to $\triangleleft
  A_0$ can be defined as
$$
P(\triangleleft A_0 ) = \{ g\in G \: ; \: \sup_{n\in \N } |\mathbf{a}^ng\mathbf{a}^{-n}| < + \infty \,
,\, \forall \mathbf{a}\in \triangleleft A_0 \}\, ,
$$ and the unipotent group corresponding to $\triangleleft A_0$,
$$
U(\triangleleft A_0 ) = \{ g\in G \: ; \: \lim_{n\to \infty } \mathbf{a}^ng\mathbf{a}^{-n} =e \, ,\,
\forall \mathbf{a}\mbox{ in the relative interior of } \triangleleft A_0 \}\, .
$$
\me

\Notat If $\triangleleft A_0^{op}$ is the opposite wall, we denote $U(\triangleleft A_0^{op} )$ by
$U_+(\triangleleft A_0)$.

\me

We have that $P(\triangleleft A_0 )=C_G(A_0)U(\triangleleft A_0 )=U(\triangleleft A_0 )C_G(A_0)$,
$U(\triangleleft A_0 )$ is the unipotent radical of $P(\triangleleft A_0 )$, and $P(\triangleleft A_0
)$ is the normalizer of $U(\triangleleft A_0 )$ in $G$.

Now let ${\mathcal{A}}=(\mathbf{a}_t)$ be a one-parameter subgroup of $G$ composed of hyperbolic
elements and let ${\mathcal{A}}^{+}$ be the positive sub-semigroup $(\mathbf{a}_t)_{t\geq 0}$. Let $r$
be a geodesic ray in $X$ such that $r(t)=r(0)\mathbf{a}_t$ for every $t\geq 0$. We consider $A_0$
either the minimal singular torus containing ${\mathcal{A}}$ or, if no such torus exists, the unique
maximal
 torus containing ${\mathcal{A}}$. We have the equality $C_G({\mathcal{A}})=C_G(A_0)$.
  If $A_0$ has dimension
one we call the one-parameter group ${\mathcal{A}}$, the semigroup ${\mathcal{A}}^{+}$ and the geodesic
ray $r$ \emph{maximal singular}.

 Let $\triangleleft A_0$ be the wall/Weyl chamber
  containing ${\mathcal{A}}^{+}\setminus \{ e \}$ in its relative
  interior. We denote $P(\triangleleft A_0 )$, $U(\triangleleft A_0 )$
  and $U_+(\triangleleft A_0 )$ also by $P(r)$, $U(r)$ and $U_+(r)$,
  respectively, and we call them the {\it{parabolic}}, the {\it{unipotent}} and the
\textit{opposite (expanding) unipotent group} of the ray $r$.

\me

Another way of defining the parabolic subgroups is as follows:

$$
P(r)= \{g\in G \; ;\; rg \sim r\}\, , \; P^0(r) = \{g\in P(r)\; ;\; r(0)g\in H(r)\}\, .
$$

The latter definition justifies calling $P^0(r)$ the horospherical group of $r$.

For instance, the parabolic group of $r_i$ is the group
\begin{displaymath}
P(r_i)=\left\{ \left(
\begin{array}{cc}
M_1 & 0 \\
N & M_2
\end{array}
\right)\in \slsr \; ;\;  M_1\in GL(s-i,\R ),\; M_2\in GL(i,\R ),\; N\in M_{i\times (s-i)}(\R )
\right\}\; .
\end{displaymath}

The horospherical subgroup is

\begin{displaymath}
P^0(r_i)=\left\{ \left(
\begin{array}{cc}
\epsilon M_1 & 0 \\
N & \epsilon M_2
\end{array}
\right) \; ;\;  M_1\in SL(s-i,\R ),\; M_2\in SL(i,\R ),\; \epsilon \in \{ \pm 1\},\; N\in M_{i\times
(s-i)}(\R ) \right\}\; .
\end{displaymath}

The opposite unipotent group is
\begin{equation}\label{uri}
U_+(r_i)=\left\{ \left(
\begin{array}{cc}
 Id_{s-i} & N \\
0 &  Id_i
\end{array}
\right)
 \; ;\; N\in M_{(s-i)\times i}(\R ) \right\}\; .
\end{equation}

\subsection{Boundary at infinity and Busemann functions of $\pps$}\label{section:bdary}

If $W$ is a Weyl chamber or a wall in $\pps$ then its boundary at infinity $W(\infty )$ is a spherical
simplex in $\partial_\infty \pps$, also called \textit{spherical chamber} or respectively
\textit{spherical wall}. These simplices cover $\partial_\infty \pps$ and determine a structure of
spherical building on it (\cite[Chapters 15,16]{Mostow:rigidity}, \cite[Appendix 5]{BGS}).

The group $\slsr$ acts on $\partial_\infty \pps$ on the right by automorphisms of spherical building.
In fact it coincides with the group of automorphisms of the spherical building $\partial_\infty \pps$.
The fundamental domain of the action of $\slsr$ on $\partial_\infty \pps$ is $W_0(\infty )$, hence one
can define a projection $\mathrm{sl}:
\partial_\infty X \to W_0 (\infty )$. The image $\mathrm{sl} (\xi )$ of every point
$\xi $ in $\partial_\infty X$ is called \textit{the slope of }$\xi$. The \textit{slope of a geodesic
ray} $r$ is the slope of $r(\infty )$.

Given a point $\xi$ in the relative interior of a spherical wall $W(\infty )$, where $W=x \triangleleft
A_0 $, the stabilizer of $\xi$ is the stabilizer of the whole wall $W(\infty )$, and it is the
parabolic group $P(\triangleleft A_0)$.

\me

The boundary at infinity $\partial_\infty {\mathcal{P}}_s$ can in fact be identified to the geometric
realization of the spherical building of flags in $\R^s$. Indeed, the complex of incidence of the flags
in $\R^s$ is a spherical building according to \cite{Tits:BN}; according to
\cite{KleinerLeeb:buildings} it can be realized geometrically as a CAT(1)-spherical complex. The sense
of the above statement is that $\partial_\infty {\mathcal{P}}_s$ endowed with the Tits metric, a
definition of which can be found in \cite{BGS}, is isomorphic and isometric to the geometric
realization of the spherical building of flags. Via this identification, the statement that $\slsr$
coincides with the group of automorphisms of the spherical building $\partial_\infty \pps$ becomes the
Fundamental Theorem of Projective Geometry.

Also via the above identification, $r_1(\infty )=\langle e_s \rangle $ and more generally $r_i(\infty
)$ is the subspace $\langle e_{s-i+1},\dots ,e_s\rangle$, for $i\in \{ 1,2,\dots s-1 \}$. The spherical
chamber $W_0 (\infty )$ is identified to the flag $\langle e_s \rangle \subset \cdots \subset \langle
e_{s-i+1},\dots ,e_s\rangle \subset \cdots \subset \langle e_{2},\dots ,e_s\rangle$.

A maximal singular ray $r$ has slope $r_i(\infty )$ if and only if $r(\infty )$ is a linear subspace of
dimension~$i$.

Given a flag $\mathcal{F}: V_1\subset \dots \subset V_k$ in $\R^s$ and a matrix $B\in GL(s,\R )$ we
denote by $B \mathcal{F}$ the flag $B (V_1)\subset \dots \subset B (V_k)$.

\begin{remark}\label{action}
The isometric action to the right $\Phi$ of $\slsr$ on $\calp_s$ induces the action to the right $\Phi
$ on $\partial_\infty {\mathcal{P}}_s$ identified to the spherical building of flags in $\R^s$, defined
by $\Phi (B , \mathcal{F} ) = B\iv \mathcal{F}$, where $\mathcal{F}$ is an arbitrary flag.
\end{remark}

The Busemann functions of $\pps$ have been computed in \cite[$\S 3.2$]{Drutu:quadrics}. We recall here
some of the results.

\begin{lemma}\label{Bus}
Let $Q$ be a positive definite quadratic form of determinant 1 on $\R^s$, let $Q_i$ be its restriction
to $\langle e_{s-i+1},\dots ,e_s\rangle $ and let $\det Q_i$ be the determinant of $Q_i$ in the basis
$\{ e_{s-i+1},\dots ,e_s\} $. Then
$$
f_{r_i}(Q)=\sqrt{\frac{s}{(s-i)i}} \ln \det Q_i \; .
$$
\end{lemma}

\bi

In particular
$$
f_{r_1}(Q)=\sqrt{\frac{s}{s-1}}\, \ln Q(e_s) \; \mbox{ and }\; f_{r_{s-1}}(Q)=\sqrt{\frac{s}{s-1}}\,
\ln Q^*(e_1)\; ,
$$
where $Q^*$ is the ``dual quadratic form'', that is the quadratic form whose matrix in the canonical
basis is $M_Q^*$, if $M_Q$ is the matrix of $Q$.

\begin{lemma}\label{bvector}
Let $d$ be a line in $\R^s$ and let $v$ be a non-zero vector on $d$.
\begin{itemize}
\item[(i)] The function $f_v :
\mathcal{P}_s \to \R $, defined by $ f_v(Q)=\sqrt{\frac{s}{s-1}}\, \ln Q(v)\, , $ is a Busemann
function of basepoint $d$.
\item[(ii)] Every Busemann function of basepoint $d$ is of the form $f_w$, where $w\in d$, $w\neq 0$.
\end{itemize}
\end{lemma}

A similar argument gives the following.

\begin{lemma}\label{bhiperpl}
Let $\hh$ be a linear hyperplane in $\R^s$ and let $v$ be a non-zero vector orthogonal to~it.
\begin{itemize}
\item[(i)] The function $f_v^* :
\mathcal{P}_s \to \R $ defined by $ f_v^*(Q)=\sqrt{\frac{s}{s-1}}\, \ln Q^*(v)\, , $ is a Busemann
function of basepoint $\hh$.
\item[(ii)] Every Busemann function of basepoint $\hh$ is of the form $f_w^*$, where $w\neq 0$ is orthogonal to $\hh$.
\end{itemize}
\end{lemma}

\me

We have that $f_{r_1}=f_{e_s}$ and $f_{r_{s-1}}=f^*_{e_1}$.

\subsection{The locally symmetric space  ${\mathcal{P}}_{s}/ \slsz$}\label{lssp}

The subgroup $\Gamma =\slsz$ is an irreducible lattice in $\slsr$, therefore $\slsr /\Gamma$ has a
$\slsz$-invariant finite measure. The $\Q$-rank $\mathbf{r}$  of $\Gamma$ is the same as the $\R$-rank
of $\slsr$, that is $s-1$. The quotient space ${\mathcal{P}}_{s}/ \Gamma$ is a non-compact locally
symmetric space of finite volume.

\me

\noindent \textit{Notations}: We denote ${\mathcal{P}}_{s}/ \Gamma$ by $\mathcal{T}_s$; we denote by
$\pr $ the projection of ${\mathcal{P}}_{s}$ onto $\mathcal{T}_s$, and by $\pr_G$ the projection of
$\slsr$ onto $\slsr /\Gamma$.

\me

Note that the space $\slsr /\Gamma$ can be naturally identified to the space of lattices of covolume
$1$ in $\R^s$, via the map $B\cdot \Gamma \mapsto B\cdot \Z^s$, where $B\in \slsr$. Consequently the
space $\mathcal{T}_s =SO(s)\, \backslash \slsr /\Gamma$ can be identified to the space of lattices of
covolume $1$ in $\R^s$ up to solid rotations preserving orientation. On the other hand, since
$\mathcal{T}_s={\mathcal{P}}_{s}/ \Gamma$, this quotient space can also be seen as the space of
positive definite quadratic forms of determinant $1$ up to the equivalence relation $Q_1 \simeq Q_2
\Leftrightarrow Q_1=Q_2 \circ B$ for some $B\in \slsz$.

\me

The projection $\pr$ restricted to the Weyl chamber $W_0$ is an isometry. Therefore $\overline{W}_0=\pr
(W_0)$ is an isometric copy of $W_0$ in $\tts$. Moreover, $\mathcal{T}_s$ is at finite Hausdorff
distance from $\overline{W}_0$. For details see \cite{Siegel:quotient} and \cite{Borel:IntroArith}.

We denote by $\bar{r}$ the projection of a ray $r$ in $W_0$.

Given a geodesic ray $r$ contained into $W_0$, the height into the end of $\tts$ can be measured by the
Busemann function $f_{\bar{r}}$ of $\bar{r}$. Moreover the following holds:

\begin{lemma}[\cite{Drutu:quadrics}, Remark 2.5.1, (1)]\label{horob}
For $a<0$ with $|a|$ large enough, the projection of the horoball $\pr (Hb_a(r))$ is the horoball
$Hb_a(\bar{r} )$.
\end{lemma}

\me

This and Lemma \ref{bvector} imply that for $a<0$ with $|a|$ large enough the projection of
$Hb_{e_s}^a$ is $Hb_a (\bar{r}_1)$ and its pre-image is $\bigcup_{v\in \pri^s} Hb_v^a$. Likewise
$Hb_{e_1^*}^a$ projects onto $Hb_a (\bar{r}_{s-1})$ and its pre-image is $\bigcup_{v\in \pri^s}
Hb_{v^*}^a$.

\me

\noindent \textit{Notation}: For simplicity, we denote in what follows the Busemann function $f_{r_i}$
on $\pps$ by $f_i$, and the Busemann function $f_{\bar{r}_i}$ on $\tts$ by $\bar{f}_i$.

\me

According to Lemma \ref{bvector}, if $\mathcal{T}_s$ is seen as the space of lattices of covolume $1$
in $\R^s$ up to solid rotations then $\bar{f}_1$ is the function associating to every lattice $\Lambda$
in $\R^s$ the value $\sqrt{\frac{s}{s-1}}\, \ln \|w \|_e$, where $w$ is a shortest non-zero vector in
$\Lambda$ with respect to the Euclidean norm $\|\cdot \|_e$.

If $\mathcal{T}_s$ is seen as the space of positive definite quadratic forms of determinant $1$ up to
the equivalence relation $\simeq$ then $\bar{f}_1$ is the function associating to every equivalence
class of quadratic forms $[Q]$ the value $\sqrt{\frac{s}{s-1}}\, \ln \lambda_1(Q)$, where
$\lambda_1(Q)$ is the first minimum of $Q$ with respect to the lattice $\Z^s$.

Likewise, using Lemma \ref{bhiperpl}, the function $\bar{f}_{s-1}$ can be seen either as the function
associating to every lattice $\Lambda$ the value $\sqrt{\frac{s}{s-1}}\, \ln vol^*_e$, where $vol^*_e$
is the minimal Euclidean covolume of a subgroup of $\Lambda$ which is a lattice in a hyperplane of
$\R^s$; or as the function associating to every equivalence class of quadratic forms $[Q]$ the value
$\sqrt{\frac{s}{s-1}}\, \ln vol^*_Q$, where $vol^*_Q$ is the minimal covolume with respect to $Q$ of a
subgroup of $\Z^s$ which is a lattice in a hyperplane of $\R^s$.

\me

The following result turns out to be a geometric version of a Transference principle.

\begin{proposition}[comparison of Busemann functions on $\tts$]\label{geomtransf}
The following inequality holds on $\tts$:
\begin{equation}\label{f1s}
(s-1)\bfs -O(1) \leq \bf1 \leq \frac{1}{s-1}\bfs +O(1)\, .
\end{equation}
\end{proposition}

\begin{remark}
The statement in Proposition \ref{horob} can be easily seen on a picture. See for instance Figure
\ref{quotient} where the case $s=3$ is represented. Note that in this case the Weyl chamber $W_0$ is
known to be an Euclidean sector of angle $\frac{\pi }{3}$.
\end{remark}

\proof The fact that $\tts$ is at finite Hausdorff distance from $\overline{W}_0$ and Lemma \ref{blip}
imply that it suffices to prove inequality (\ref{f1s}) for the restrictions of $\bfs$ and $\bf1$ to
$\overline{W}_0$. Or $\overline{W}_0$ can be identified to the following polytopic cone in $\R^{s}$:
$$
\{ (t_1,t_2,\dots , t_s)\; ;\; t_1+t_2+\cdots +t_s=0,\; t_1\geq t_2\geq\cdots \geq t_s \}\, .
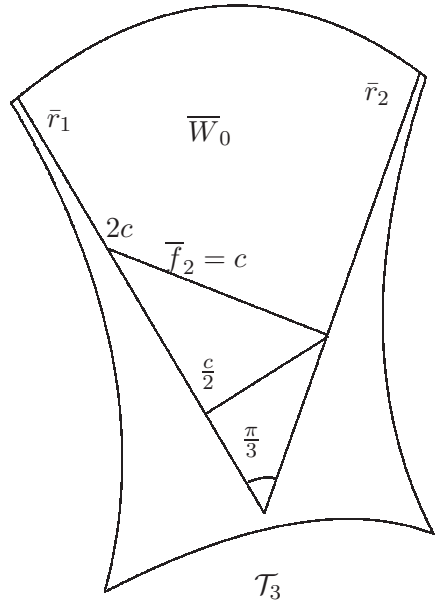
$$

\unitlength 1mm 
\linethickness{0.4pt}
\ifx\plotpoint\undefined\newsavebox{\plotpoint}\fi 
\begin{picture}(84.25,96.5)(0,0)
\put(84.25,75.25){\line(0,1){0}}
\multiput(25.5,73)(.0337281153,-.056900103){971}{\line(0,-1){.056900103}}
\multiput(58.25,17.75)(.033739837,.095528455){615}{\line(0,1){.095528455}}
\put(79,76.5){\line(0,1){0}}
\qbezier(24.5,72.5)(45.25,37.63)(37,7.25) \qbezier(79.75,75.75)(67.75,39.63)(80.75,15)
\qbezier(80.75,15)(64.13,21.38)(37,7.25) \qbezier(24.5,72.25)(52.88,96.5)(79.75,75.75)
\qbezier(56.25,22)(59,22.88)(59.75,22.25)
\multiput(37.25,53)(.085777126,-.03372434){341}{\line(1,0){.085777126}}
\put(31,70){\makebox(0,0)[cc]{$\bar{r}_1$}} \put(73.25,73.5){\makebox(0,0)[cc]{$\bar{r}_2$}}
\put(56.5,27.25){\makebox(0,0)[cc]{$\frac{\pi}{3}$}}
\multiput(66.75,41.25)(-.053453947,-.033717105){304}{\line(-1,0){.053453947}}
\put(50.5,31){\line(0,1){0}}
\put(50.5,51.75){\makebox(0,0)[cc]{$\overline{f}_2=c$}} \put(39,55.75){\makebox(0,0)[cc]{$2c$}}
\put(50.75,36.25){\makebox(0,0)[cc]{$\frac{c}{2}$}} \put(51,68.75){\makebox(0,0)[cc]{$\overline{W}_0$}}
\put(58.75,8){\makebox(0,0)[cc]{$\tt_3$}}
\end{picture}

\begin{figure}[!ht]
\centering \caption{Inequality (\ref{f1s}) in case $s=3$.} \label{quotient}
\end{figure}

According to Example \ref{brn}, for any $i\in \{ 1,2,\dots s-1 \}$ the Busemann function $\bar{f}_i$
restricted to $\overline{W}_0$ coincides, via this identification, with $-\langle \cdot | v_i \rangle$,
where
\begin{equation}\label{vi}
v_i=(\underbrace{\lambda_i,\dots ,\lambda_i}_{s-i\mbox{ times}}, \underbrace{-\mu_i,\dots , -\mu_i
t}_{i\mbox{ times}})\mbox{, with }\lambda_i=\sqrt{\frac{i}{s(s-i)}}\mbox{ and
}\mu_i=\sqrt{\frac{s-i}{si}}\, .
\end{equation}

Any two rays in $\overline{W}_0$ with same vertex as $\overline{W}_0$ form an angle strictly smaller
than $\frac{\pi }{2}$. This can be verified in this case by direct computation, and it also follows
from general results stating that in a Weyl chamber any two rays with origin in its vertex form an
angle smaller or equal to $\frac{\pi }{2}$, and the equality case may appear if and only if the
corresponding symmetric space is reducible, i.e. it decomposes as a cartesian product. See for instance
\cite{KleinerLeeb:buildings} where the latter result is explained in a nice and geometric way.

It follows that any horosphere $H_a(\bar{r})$ with $\bar{r}$ a ray in $\overline{W}_0$ intersects
$\overline{W}_0$ in a finite polytope. In particular it is the case for a horosphere defined by $\bfs =
-c$, with $c$ a large enough positive constant $c$.

The maximum and minimum of $\bf1$ on the above polytope must be attained in one of the vertices, since
$\bf1$ restricted to the polytope coincides with the linear function $-\langle \cdot | v_i \rangle $.
Or the vertices are here the intersections of the horosphere $H_{-c}(\bar{r}_{s-1})$ with $\bar{r}_{i}$
for all $\is$. One easily sees that they are $\bar{r}_{i} (t_i)$, with $t_i=\frac{c}{\langle v_i |
v_{s-1} \rangle}$, which by the above identification of $\overline{W}_0$ to a polytopic cone in $\R^s$,
become $t_iv_i$.

Now $\bf1 (\bar{r}_{i} (t_i) )= c \frac{\langle v_i | v_{1} \rangle}{\langle v_i | v_{s-1} \rangle}$.
Elementary computations give $\langle v_i | v_{1} \rangle = \sqrt{\frac{s-i}{i(s-1)}}$, and  $\langle
v_i | v_{s-1} \rangle = \sqrt{\frac{i}{(s-i)(s-1)}}$, hence $\bf1 (\bar{r}_{i} (t_i) )=c \left(
\frac{s}{i}-1 \right)$. We conclude that the maximum of $\bf1$ on $\overline{W}_0\cap
H_{-c}(\bar{r}_{s-1})$ is attained for $i=1$, and it is $s-1$, while the minimum is attained for
$i=s-1$, and it is $\frac{1}{s-1}$. This implies that the inequalities in (\ref{f1s}) hold on
$\overline{W}_0$ without the $O(1)$ terms.
\endproof

\me

\section{Diophantine approximation and excursions of geodesic rays}\label{sdioph}

\subsection{Diophantine approximation for families of forms}

In what follows we study from the Diophantine approximation viewpoint families of $\ell$ linear forms
in $m$ variables, and their transposed family of $m$ linear forms in $\ell$ variables:
$$
L_i(x_1,...,x_m) =\sum_{j=1}^m a_{ij} x_j\, ,\; \;  M_j(y_1,..., y_\ell) = \sum_{i=1}^\ell a_{ij} y_i
\, .
$$

\me

\Notat \quad We denote by $L$ the matrix $(a_{ij})_{1\leq i\leq \ell , 1\leq j\leq m}$ and by $M$ its
transposed. We also denote the sum $\ell +m$ by $s$.

\me


Given an approximating function $\phi $ we consider the set of $\phi$-approximable families of $\ell$
linear forms in $m$ variables
\begin{equation}\label{lphi}
\cl_\phi = \{ L\in M_{\ell \times m } (\R )\; ;\;  |L_i(\bar{q}) - p_i |\leq \phi \left( \| \bar{q}
\|_\ma \right) \mbox{ for infinitely many }(\bar{p}, \bar{q})\in \pri^s \}\, .
\end{equation}

Note that if $L$ satisfies the hypothesis in Theorem \ref{functgenkhintch} then $L$ is in $\cl_\phi$.

Similarly we define the set
\begin{equation}\label{mpsi}
\mm_\psi = \{M\in M_{m \times \ell } (\R )\; ;\; |M_j(\bar{a}) - b_i |\leq \psi \left( \| \bar{a}
\|_\ma \right)\mbox{ for infinitely many }(\bar{a}, \bar{b})\in \pri^s\}\, .
\end{equation}

Both the set $\cl_\phi$ and the set $\mm_\psi$ can be related to sets of geodesic rays of the same
slope as $r_m$ and rising further and further in the cusp.

\subsection{Two collections of geodesic rays}

\me

\Notat To simplify the formulas, we use the notation $\eta$ for the constant $\sqrt{\frac{s}{s-1}}$. We
shall also continue using $\lambda_i$ and $\mu_i$ to designate the constants defined in (\ref{vi}) for
$\is \, ,\, i\neq m$, while for $i=m$ we shall drop henceforth the index, and simply write $\lambda$
and $\mu$.

\me

We introduce now two sets of geodesic rays, we explain why their definition is natural, and in the end
we explain how these two sets relate to sets of type $\cl_\phi$ and respectively $\mm_\psi$.

\me

Consider a (strictly) increasing function $\varphi :[a, +\infty ) \to [b, +\infty )$ and for $k =
1,s-1$ define the following set of unipotents:
\begin{equation}\label{rr1}
\mathcal{R}_{\varphi }^k = \left\{ \mathbf{u}\in U_+(r_m) \; ;\; -\bar{f}_k \left(\pr
\left(r_m(t)\mathbf{u}\right)\right)\geq \alpha_k t - \varphi (t) \mbox{ infinitely many times as }t\to
\infty \right\}\, ,
\end{equation} where $\alpha_k = \langle v_k| v_m\rangle$, with $v_i$ the vectors defined in
(\ref{vi}). Thus
\begin{equation}\label{alphas}
\alpha_1 = \langle v_1| v_m\rangle = \sqrt{\frac{m}{\ell (s-1)}} = \eta \mu \, ,\; \; \alpha_{s-1} =
\langle v_{s-1}| v_m\rangle = \sqrt{\frac{\ell}{m (s-1)}} = \eta \lambda \, .
\end{equation}

The set $\mathcal{R}_{\varphi }^k$, though a set of unipotents, can be identified to a set of rays of
same slope as $r_m$ in $\pps$ via the bijection
$$
\bu \mapsto r_m\bu \, .
$$

\begin{figure}[!ht]
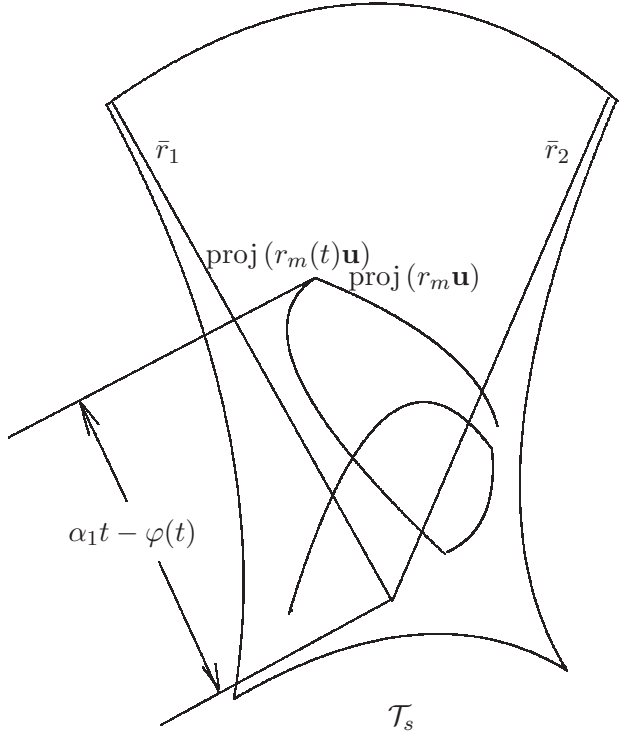

\centering \caption{The set $\mathcal{R}_{\varphi }^1$.} \label{rphi}
\end{figure}

These rays have the property that their projection in $\tts$ rises infinitely many times
 in the cusp at height at least $ \alpha_k t - \varphi (t)$, where the height is measured with respect to
 the ray $\bar{r}_k$ and $t$ is the time at which that height is attained (see Figure \ref{rphi}).

Several explanations are needed concerning the choice of defining $\mathcal{R}_{\varphi }^k$ as a set
of unipotents, and the form of the function measuring the height.

\begin{remark}[why a set of unipotents]\label{last}
All geodesic rays having the same slope as $r_m$ are in the orbit $r_m G$, and the set $ P(r_m)
U_+(r_m)$ is open Zariski dense in $\slsr$.

If a geodesic ray has a projection on $\tts$ moving away into the cusp infinitely many times with
height measured by the function $\alpha_k id -\phi $ with respect to the ray $\bar{r}_k$, then any
geodesic ray asymptotic to it has the same property, up to a bounded perturbation of the height. Thus
if $\mathbf{u}\in \mathcal{R}_{\varphi }^k$ then any geodesic ray $\rho$ in $r_m P(r_m)\mathbf{u}$ has
the property that $-\bar{f}_k\left(\pr \left(\rho(t)\right)\right)\geq \alpha_k t - \varphi (t)-C
\mbox{ infinitely many times as }t\to \infty $ for some constant $C=C(\rho )$.

We may therefore say that the set $\mathcal{R}_{\varphi }^k$ deals with all rays with same slope as
$r_m$ and ascending speed in the cups measured by the function $\alpha_k id -\phi +O(1)$ with respect
to $\bar{r}_k$, with the exception of an algebraic variety.
\end{remark}

\unitlength 1mm 
\linethickness{0.4pt}
\ifx\plotpoint\undefined\newsavebox{\plotpoint}\fi 
\begin{picture}(95,116.75)(0,0)
\multiput(27.75,90.5)(.0337104072,-.0601809955){1105}{\line(0,-1){.0601809955}}
\multiput(65,24)(.0337045721,.0785463072){853}{\line(0,1){.0785463072}}
\qbezier(27,90)(50,47.38)(44,11.25) \qbezier(94.75,90.5)(73,38.88)(88.25,14.75)
\qbezier(88.25,14.75)(69.88,26.13)(44,11) \qbezier(27,90)(64.5,116.75)(95,90.5)
\put(51.5,22.75){\line(1,0){.25}}
\qbezier(51.25,22.25)(63.75,63.75)(78.25,44.25) \qbezier(78.25,44.25)(79.25,33.63)(72.25,30.5)
\qbezier(72,30.25)(42.13,58.38)(54.75,67) \qbezier(54.75,67)(76.13,58.38)(79,47.25)
\put(32.5,6.75){\line(0,1){0}}
\multiput(54.75,67)(-.06468254,-.033730159){630}{\line(-1,0){.06468254}}
\multiput(65,24.25)(-.061871227,-.033702213){497}{\line(-1,0){.061871227}}
\multiput(23.5,50.5)(.03360215,-.07123656){186}{\line(0,-1){.07123656}}
\multiput(34,29.25)(.03361345,-.07352941){238}{\line(0,-1){.07352941}}
\put(42,11.75){\line(0,1){0}}
\multiput(23.5,50.5)(.0328947,-.1184211){38}{\line(0,-1){.1184211}}
\multiput(23.75,50.25)(.0335821,-.0410448){67}{\line(0,-1){.0410448}}
\multiput(41.5,15.75)(.033333,-.233333){15}{\line(0,-1){.233333}}
\put(42,12.25){\line(0,1){0}}
\multiput(42,12.25)(-.0335821,.0447761){67}{\line(0,1){.0447761}}
\put(39.75,15.25){\line(0,1){0}}
\put(35.25,83.75){\makebox(0,0)[cc]{$\bar{r}_1$}} \put(87,83.75){\makebox(0,0)[cc]{$\bar{r}_2$}}
\put(68,67){\makebox(0,0)[cc]{$\pr \left(r_m\mathbf{u}\right)$}}
\put(30.5,33){\makebox(0,0)[cc]{$\alpha_1 t -\varphi (t)$}} \put(51.25,70){\makebox(0,0)[cc]{$\pr
\left(r_m(t)\mathbf{u}\right)$}} \put(66.25,8.5){\makebox(0,0)[cc]{$\tts$}}
\end{picture}

\begin{lemma}[why the ascending function $\alpha_k id -\phi $]\label{lemdiv}
If a geodesic ray $\rho$ in $\pps$ has the property that $-\bar{f}_k\left(\pr
\left(\rho(t)\right)\right)\geq \alpha_k t +O(1) \mbox{ infinitely many times as }t\to \infty$ then
$-\bar{f}_k\left(\pr \left(\rho(t)\right)\right)= \alpha_k t +O(1)$, and $\rho$ is asymptotic to a ray
$\rho'$ contained in the same Weyl chamber as a lift of $\bar{r}_k$ (i.e. a ray in $r_k \Gamma$).
\end{lemma}

Lemma \ref{lemdiv} thus implies that $\alpha_k\, id$ is the maximal possible ascending function with
respect to $\bar{r}_k$.

\proof According to the argument in Remark \ref{last}, the case when $\rho$ is in $r_mP(r_m)U_+(r_m)$
can be reduced to the case when $\rho$ is in $r_m U_+(r_m)$.

The other cases can likewise be reduced to sets of rays of the form $r_m U_+(r_m)w$, where $w$ is one
of the elements of the Weyl group corresponding to $F_0$ (that is, and element in $Z_G(A_0)\backslash
N_G(A_0)$, hence an element that can be realized as a permutation matrix). But applying a permutation
matrix does not change the data in an important way. Thus, the arguments that we give below for a ray
in $r_m U_+(r_m)$ also works for rays in $r_m U_+(r_m)w$. We leave this as an exercise to the reader.

\me

\emph{Case $k=1$}. \quad Let $\bu$ be a unipotent in $U_+(r_m)$. The fact that $-\bf1\left(\pr
\left(r_m(t)\bu \right)\right)\geq \alpha_1 t +O(1) \mbox{ infinitely many times as }t\to \infty$
implies that for infinitely many $(\bar{p},\bar{q})\in \pri^s$ the value $f_{(\bar{p},\bar{q})}
(r_m(t)\bu )$ is at most $-\alpha_1 t +O(1)$. According to Lemma \ref{bvector}, (i), this is equivalent
to the fact that $r_1(t)\bu (\bar{p},\bar{q})$ is at most $C e^{-\frac{\alpha_1 t}{\eta }}= C e^{-\mu
t}$, for some constant $C$.

Assume that $\bu$ is as in the formula (\ref{uri}), for some matrix $N\in M_{\ell \times m}(\R)$. Then
the above implies that
$$
e^{\lambda t} \|\bar{p} + N \bar{q}\|^2 + e^{-\mu t} \|\bar{q}\|^2 \leq C e^{-\mu t} \Rightarrow
\|\bar{p} + N \bar{q}\|\ll e^{-\frac{\lambda +\mu }{2}t}\, ,\; \; \|\bar{q}\| \ll 1 \, .
$$

The last inequality implies that up to taking a subsequence we may assume that $\bar{q}$ is constant.
Then, as $\bar{p} + N \bar{q}$ now varies in the lattice $\Z^\ell + N \bar{q}$, the first inequality,
when $t$ is large enough, implies that $\bar{p} + N \bar{q} =0$, in particular $\bar{p}$ is a constant.
For the fixed vector $(\bar{p},\bar{q})\in \pri^s$ we have that $f_{(\bar{p},\bar{q})} (r_m(t)\bu )=
\eta \ln [e^{\lambda t} \|\bar{p} + N \bar{q}\|^2 + e^{-\mu t} \|\bar{q}\|^2]= -\eta \mu t + 2\eta \ln
\|\bar{q}\|$, hence $-\bf1 \left(\pr \left(r_m(t)\bu \right)\right)= \alpha_1 t +O(1)$.

From the equality $\bar{p} + N \bar{q} =0$ also follows that
$$
\left(
\begin{array}{cc}
Id_{s-i} & N \\
0 & Id_i
\end{array}
\right)\cdot \left(
\begin{array}{c}
\bar{q} \\
\bar{p}
\end{array}
\right) = \left(\begin{array}{c}
0 \\
\bar{p}
\end{array}
\right)\, .
$$

This can be rewritten as $\bu (\bar{q}, \bar{p}) = (0, \bar{p}) $.

Note that, since the vector $(0, \bar{p})$ is in the subspace $\langle e_{\ell+1}, ..., e_s\rangle$,
which can be identified to $r_m(\infty )$ (see Section \ref{section:bdary}), the vector $(0, \bar{p})$
can be written as $\rho (\infty )$ for some ray in the same Weyl chamber as $r_m$.

Take $\gamma \in \Gamma$ such that $\gamma\iv (e_s) = (\bar{q}, \bar{p})$. This is the same as writing
that $r_1 (\infty )\gamma $ is $(\bar{q}, \bar{p})$. Therefore $r_1 (\infty )\gamma \bu\iv $ is $(0,
\bar{p})$. This implies that $r_1 (\infty )\gamma = \rho (\infty ) \bu $. Now $\rho \bu$ is in the same
Weyl chamber as $r_m\bu$, and $r_1 \gamma$ projects in $\tts$ onto $\bar{r}_1$, hence $\rho\bu$
projects to a geodesic ray asymptotic to $\bar{r}_1$.

\me

\emph{Case $k=s-1$}. \quad Assume now that $\bu \in U_+(r_m)$ is such that $-\bfs \left(\pr
\left(r_m(t)\bu \right)\right)\geq \alpha_{s-1} t +O(1) \mbox{ infinitely many times as }t\to \infty$.
Then there exist infinitely many $(\bar{a},\bar{b})\in \pri^s$ such that $f^*_{(\bar{a},\bar{b})}
(r_m(t)\bu )$ is at most $-\alpha_{s-1} t +O(1)$. Equivalently (with the same form of $\bu$ as in Case
1) $e^{-\lambda t} \|\bar{a}\|^2_e + e^{\mu t} \|N^T \bar{a} +\bar{b}\|^2 \leq C
e^{-\frac{\alpha_{s-1}t}{\eta }}= e^{-\lambda t }$.

It follows that eventually by taking a subsequence we may assume that $\bar{a}$ is fixed. Also, since
$N^T \bar{a} +\bar{b} $ is in the lattice $ N^T \bar{a} +\Z^m$ and since $\|N^T \bar{a} +\bar{b}\| \ll
e^{-\frac{\lambda +\mu }{2}t}$ with $t\to \infty $ it follows that $N^T \bar{a} +\bar{b}=0$ for some
$\bar{b} \in \Z^m$. For the fixed vector $(\bar{a},\bar{b})\in \pri^s$ thus found, a straightforward
application of Lemma \ref{bhiperpl} gives that $f^*_{(\bar{a},\bar{b})} (r_m(t)\bu ) = -\bfs \left(\pr
\left(r_m(t)\bu \right)\right) = -\alpha_{s-1} t +O(1)$.

Thus, $\bu\iv$ applied to the rational hyperplane of coefficients $(-\bar{a},\bar{b})$ gives the
rational hyperplane of coefficients $(-\bar{a}, 0)$. The subspace  $\langle e_{\ell+1}, ...,
e_s\rangle$ is contained in this hyperplane, hence the hyperplane of normal vector $(-\bar{a}, 0)$ is
$\rho (\infty )$ for some $\rho$ in the same Weyl chamber as $r_m$ (and of same slope as $r_{s-1}$).

Take $\gamma \in \Gamma$ such that $\gamma^T (e_1)= (-\bar{a},\bar{b})$, equivalently $r_{s-1}(\infty
)\gamma = (-\bar{a},\bar{b})$. Then $r_{s-1}(\infty )\gamma \bu\iv = \rho (\infty )$. Hence $\rho \bu$
projects to a geodesic ray asymptotic to $\bar{r}_{s-1}$, and it is in the same Weyl chamber as
$r_m\bu$.\endproof

Proposition \ref{geomtransf} immediately implies the following relation between the two sets of
ascending rays defined in (\ref{rr1}).

\begin{lemma}\label{relr}
\begin{enumerate}
    \item Let  $\varphi :[a, +\infty ) \to [b, +\infty )$ be a strictly increasing function. Then
    \begin{equation}\label{1ins}
   \calr_\varphi^1 \subset \calr_\psi^{s-1} \mbox{ with }\psi (t) = \eta \left( \lambda - \frac{\mu }{s-1}
   \right)t + \frac{1}{s-1}\varphi (t) +O(1)\, .
\end{equation}

    \item  Let  $\psi :[a', +\infty ) \to [b', +\infty )$ be a strictly increasing function. Then
    \begin{equation}\label{sin1}
   \calr_\psi^{s-1} \subset \calr_\varphi^1 \mbox{ with }\varphi(t) = \eta \left( \mu - \frac{\lambda }{s-1}
   \right)t + \frac{1}{s-1}\psi (t) +O(1)\, .
\end{equation}
\end{enumerate}
\end{lemma}

\subsection{Relation between sets of linear forms and sets of rays}

The set of matrices $ M_{\ell \times m} (\R )$ can be naturally identified to $U_+(r_m)$, with the map
\begin{equation}\label{lmap}
L\mapsto \left(
\begin{array}{cc}
Id_{\ell} & L \\
0 & Id_m
\end{array}
\right)\, .
\end{equation}

With this identification, the set $\cl_\phi$ defined in (\ref{lphi}) becomes a subset $\tl_\phi $ in
$U_+(r_m)$.

Likewise the set of matrices $ M_{m \times \ell} (\R )$ can be naturally identified to $U_+(r_m)$ by
means of the map
\begin{equation}\label{mmap}
M\mapsto \left(
\begin{array}{cc}
Id_{\ell} & M^T \\
0 & Id_m
\end{array}
\right)\, .
\end{equation}

As previously, with this map, we identify the set $\mm_\psi$ from (\ref{mpsi}) with a subset $\tm_\psi$
in $U_+(r_m)$.

We establish the following relations between sets of well approximable linear forms and geodesic rays
rising in the cusp.

\begin{proposition}\label{lr}
 Let $\varphi :[a,\infty ) \to [b, \infty)$ be a function such that $\varphi$ and $\eta (\lambda + \mu ) id
    -\varphi$ are (strictly) increasing. Then
    \begin{equation}\label{fr1}
    \tl_{\Phi_1} \subset \calr_\varphi^1 \subset \tl_{\Phi_2}\mbox{ and }\tm_{\Phi_1} \subset \calr_\varphi^{s-1} \subset \tm_{\Phi_2}
\end{equation} where $\Phi_1(x) = \frac{1}{\sqrt{s}}x e^{-\frac{\lambda + \mu }{2}\varphi\iv (2\eta \ln
(\sqrt{2s}x))}$ and  $\Phi_2(x) = \sqrt{s}x e^{-\frac{\lambda + \mu }{2}\varphi\iv (2\eta \ln (x))}$.
\end{proposition}

\proof \emph{Case} $k=1$. \quad We prove the first inclusion. Assume that $L\in \cl_{\Phi_1}$, that is
there exist infinitely many $(\bar{p}, \bar{q})\in \pri^s$ such that $\|\bar{p} + L\bar{q}\|_\ma \leq
\Phi_1 (\|q\|_\ma)$. For each such primitive vector consider $t= \varphi\iv (2\eta \ln (\sqrt{2}
\|\bar{q}\|_e))$. Then $e^{-\mu t} \|\bar{q}\|^2_e = \frac{1}{2} e^{\frac{\varphi (t) -\alpha_1
t}{\eta}}$.

On the other hand $e^{\lambda t} \|\bar{p} + L \bar{q}\|_e^2\leq s e^{\lambda t} \|\bar{p} + L
\bar{q}\|_\ma^2 \leq e^{\lambda t} \|\bar{q}\|^2_\ma e^{-(\lambda +\mu ) t} \leq e^{\lambda t}
\|\bar{q}\|^2_e e^{-(\lambda +\mu )t}= e^{-\mu t} \frac{1}{2}e^{\frac{\varphi (t)}{\eta }}=
\frac{1}{2}e^{\frac{\varphi (t) - \alpha_1 t}{\eta}}$.

On the whole we obtain that $e^{\lambda t} \|\bar{p} + L \bar{q}\|_e^2 + e^{-\mu t} \|\bar{q}\|^2_e
\leq e^{\frac{\varphi (t) - \alpha_1 t}{\eta}}$, whence $f_{(\bar{p}, \bar{q})} (r_m(t) \bu )\leq
\varphi (t) - \alpha_1 t$.

\me

Now we prove the second inclusion. Take a unipotent $\bu $ corresponding to a matrix $L$ such that
$-\bf1 (\pr (r_m(t)\bu )) \geq \alpha_1 t-\varphi (t)$ infinitely many times as $t$ goes to infinity.

Then for infinitely many $(\bar{p}, \bar{q})\in \pri^s$ we have for some $t>0$ that
$$
e^{\lambda t} \|\bar{p} + L \bar{q}\|_e^2 + e^{-\mu t} \|\bar{q}\|^2_e \leq e^{\frac{\varphi (t) -
\alpha_1 t}{\eta}}\, .
$$

It follows that $\|\bar{q}\|^2_e \leq e^{\frac{\varphi (t)}{\eta}}$. This and the fact that $\varphi$
is increasing imply that $t\geq \varphi\iv (2\eta \ln \|\bar{q}\|_e)$.

Then $\|\bar{p} + L \bar{q}\|_e^2 \leq e^{\frac{\varphi (t)-(\lambda + \mu )t}{\eta}}$. The hypothesis
that $\eta (\lambda + \mu ) id
    -\varphi$ is increasing implies that the latter term is smaller than
    $\|\bar{q}\|_e^2 e^{-(\lambda + \mu )\varphi\iv (2\eta \ln \|\bar{q}\|_e)}$. Whence $\|\bar{p} + L \bar{q}\|_\ma \leq
    \Phi_2(\|\bar{q}\|_\ma)$.

\me

\emph{Case} $k=s-1$. \quad Take $M\in \mm_{\Phi_1}$. Then for infinitely many $(\bar{a}, \bar{b})\in
\pri^s$ we have that $\| M(\bar{a}) - \bar{b}\|_\ma \leq \Phi_1 \left( \| \bar{a} \|_\ma \right)$.

 For every  $(\bar{a}, \bar{b})$ as above let $t= \varphi\iv (2\eta \ln
(\sqrt{2} \|\bar{a}\|_e))$, equivalently such that $\|\bar{a}\|_e^2 = \frac{1}{2}e^{\frac{\varphi
(t)}{\eta }}$.

By hypothesis $\|M\bar{a} +\bar{b}\|_e^2 \leq s\|M\bar{a} +\bar{b}\|_\ma^2 \leq  \|\bar{a}\|_\ma^2
e^{-(\lambda +\mu ) t} \leq \frac{1}{2}e^{\frac{\varphi (t)}{\eta }}e^{-(\lambda +\mu ) t}\, .$

We conclude that $e^{-\lambda t}\|\bar{a}\|_e^2 + e^{\mu t} \|M\bar{a} +\bar{b}\|_e^2 \leq
e^{\frac{\varphi (t)}{\eta } - \lambda t} = e^{\frac{\varphi (t) - \alpha_{s-1}t}{\eta }}\, .$ This and
Lemma \ref{bhiperpl} imply that, if $\bu$ is the unipotent corresponding to the matrix $M$, then we may
write $f^*_{(\bar{a}, \bar{b})} (r_m(t)\bu )\leq \varphi (t) - \alpha_{s-1}t$.

For the second inclusion, assume that $\bu $ is such that $-\bfs (\pr (r_m(t)\bu )) \geq \alpha_{s-1}
t-\varphi (t)$ infinitely many times as $t$ goes to infinity.

Then there exist infinitely many $(\bar{a}, \bar{b})\in \pri^s$, and $t>0$, such that
$$
e^{-\lambda t}\|\bar{a}\|^2_e  + e^{\mu t} \|M\bar{a} +\bar{b}\|_e^2 \leq e^{\frac{\varphi (t) -
\alpha_{s-1} t}{\eta}}\, .
$$

Then $\|\bar{a}\|^2_e \leq e^{\frac{\varphi (t)}{\eta}}$, which implies that $t\geq \varphi\iv (2\eta
\ln \|\bar{a}\|_e)$.

It follows that $ \|M\bar{a} +\bar{b}\|_e^2 \leq e^{\frac{\varphi (t)-(\lambda + \mu )t}{\eta}}$. Since
$\eta (\lambda + \mu ) id
    -\varphi$ is increasing we may bound the last term by
    $\|\bar{a}\|_e^2 e^{-(\lambda + \mu )\varphi\iv (2\eta \ln \|\bar{a}\|_e)}$, and conclude that $\|M\bar{a} +\bar{b}\|_\ma \leq
    \Phi_2(\|\bar{a}\|_\ma)$.\endproof

\begin{remark}
The conditions on the function $\varphi$ are not so restrictive, in the sense that one does not really
exclude from discussion some of the ascending rays. Indeed, if a ray $\rho$ is in a set
$\calr^k_\varphi$ with a positive function $\varphi$ such that $\lim_{t\to \infty }\varphi (t)= \infty
$ (otherwise we are in the case of Lemma \ref{lemdiv}) and such that $\lim_{t\to \infty } \alpha_k t -
\varphi (t) = \infty$ (otherwise the set of rays is uninteresting) then one can choose a sequence of
parameters $t_n \to \infty $ such that $-\bar{f}_k\left(\pr \left(r_m(t_n)\mathbf{u}\right)\right)\geq
\alpha_k t_n - \varphi (t_n)$ and such that $\varphi $ restricted to the sequence $(t_n)$ is
increasing. By replacing $\varphi$ with a piecewise affine map coinciding with $\varphi$ on $(t_n)$ one
can make both $\varphi $ and $\eta (\lambda + \mu ) id
    -\varphi$ strictly increasing.
    \end{remark}

\subsection{Transference principles deduced from geometry of the locally symmetric space
$\tts$}\label{calc}

Proposition \ref{geomtransf} can be used to obtain transference principles. Since the computations
needed to deduce Theorem \ref{functgenkhintch}, (I), are more elaborate, we shall only give the
arguments needed to deduce Theorem \ref{functgenkhintch}, (II), to give an idea of how it all works.

Take a matrix $L\in \cl_\Phi$, with $\Phi (x)= x^{-\frac{m+\alpha }{\ell }}$, let $\bu$ be the
corresponding unipotent, and let $M=L^T$. According to (\ref{fr1}), $\bu$ is in $\calr^1_\varphi$ with
$\varphi$ such that
$$
x^{-\frac{m+\alpha }{\ell }-1} = e^{-\frac{\lambda + \mu }{2}\varphi\iv (2\eta \ln (\sqrt{2s}x))}\, .
$$

From this it can be deduced by a simple computation that
$$
\varphi (t) = \frac{s}{s+\alpha } \sqrt{\frac{\ell }{m (s-1)}}\, \cdot \, t +O(1)\, .
$$

Proposition \ref{geomtransf} implies that $\bu \in \calr^{s-1}_\psi$, with $\psi$ such that $\psi (t) =
\eta \left( \lambda -\frac{\mu }{s-1} \right)t + \frac{1}{s-1} \varphi (t) + O(1)$.

This gives $\psi (t) = \frac{s[m(s+\alpha ) - m -\alpha ]}{(s+\alpha )(s-1) \sqrt{\ell m (s-1)}}\,
\cdot \, t +O(1)$.

By the last inclusion of (\ref{fr1}) we then deduce that $M$ is in $\mm_\Psi$, with $\Psi \asymp
x^{1-\frac{(s+\alpha )(s-1)}{m (s+\alpha )-m - \alpha }}$. The exponent of $x$ can be rewritten as
$\frac{\ell +\beta }{m}$ with  $\beta = \frac{\ell \alpha }{m(s -1)+ (m-1)\alpha }\, .$


\providecommand{\bysame}{\leavevmode\hbox to3em{\hrulefill}\thinspace}
\providecommand{\MR}{\relax\ifhmode\unskip\space\fi MR }
\providecommand{\MRhref}[2]{%
  \href{http://www.ams.org/mathscinet-getitem?mr=#1}{#2}
} \providecommand{\href}[2]{#2}

\end{document}